\documentclass[10pt,a4paper,twoside,bibliography=totocnumbered]{article}

\usepackage{cite}

\usepackage{xcolor}

\usepackage{graphicx}

\usepackage{amsmath}
\usepackage{amssymb}

\usepackage{geometry}

\usepackage{pdflscape}

\begin{document}

% title page
\title{{\huge Steady-state distribution and stability analysis of random differential equations with uncertainties and superpositions: Application to a predator-prey model\vspace*{-0.4cm}\\ \begin{minipage}{5cm}\centering \textcolor{black}{---}\vspace*{-0.4cm}
\end{minipage}\\}{\Large \textcolor{black}{Original Research Article}}\vspace*{0.5cm}}

\author{Wolfgang Hoegele$^{1}$\vspace*{0.4cm}\\
{\normalsize $^{1}$ Munich University of Applied Sciences HM}\\{\normalsize Department of Computer Science and Mathematics}\\{\normalsize Lothstraße 64, 80335 München, Germany}\vspace*{0.4cm}\\
{\normalsize corresponding mail: \texttt{wolfgang.hoegele@hm.edu}}\vspace*{0.4cm}\\
ORCID: 0000-0002-5303-9334 \vspace*{0.4cm}}

\date{\today}

\maketitle
\thispagestyle{empty}

\textit{This manuscript is published and can be cited by}\medskip

\textit{W. Hoegele, Steady-state distribution and stability analysis of random differential equations with uncertainties and superpositions: Application to a predator-prey model,\\ International Journal of Modelling and Simulation, 2026,\\ https://doi.org/10.1080/02286203.2026.2767580}\medskip

\textit{This manuscript appears also on ArXiv.org\\
arXiv:2603.03845 [math.DS],\\ https://doi.org/10.48550/arXiv.2603.03845}

\section*{Abstract}
We present a computational framework to investigate steady-state distributions and perform a stability analysis for random ordinary differential equations driven by parameter uncertainty. Using the nonlinear Rosenzweig-MacArthur predator-prey model as a case study, we characterize the non-trivial steady state of the system and investigate its complex distribution when the parameter probability densities are multimodal mixture models with partially overlapping or separated components. Consequently, this application incorporates both parameter uncertainty and superposition. In addition, we present the stability analysis of steady states based on the eigenvalue distribution of the system's Jacobian matrix in this stochastic regime. 
The steady-state posterior density and stability metrics are computed with a recently published Monte Carlo based numerical scheme specifically designed for random equation systems. In particular, the simplicity of this stochastic extension of dynamical systems combined with a broadly applicable computational approach is demonstrated. Numerical experiments show the emergence of multimodal steady-state distributions of the predator-prey model, and we calculate their stability regions, illustrating the method's applicability to uncertainty quantification in dynamical systems.\medskip

\textbf{Keywords:} random differential equation, steady state, stability analysis, predator-prey model, Rosenzweig-MacArthur model\medskip

\newpage
\thispagestyle{empty}

%% table of contents
\tableofcontents

\section*{About the Author}

Dr. Högele is Professor of Applied Mathematics and Computational Science at the Department of Computer Science and Mathematics at the Munich University of Applied Sciences HM, Germany. His research interests are mathematical and stochastic modeling, simulation and analysis of complex systems in applied mathematics.\medskip

\thispagestyle{empty}

\newpage

\section{Introduction}

Computing steady-state solutions of ordinary differential equations (ODEs) is an important task in order to determine the long-term states toward which dynamical systems evolve after an initial transient phase of typically high volatility. The existence and sensitivity of those steady states depend highly on the specific system under consideration and its parameters. Parameter uncertainty or incomplete knowledge, represented by probability density functions, complicates this task enormously and is described by random differential equations (RDE) \cite{Jornet2023}. In this case, one can describe the steady state only stochastically by its probability density function, which is part of uncertainty quantification of dynamical systems and is closely connected to soft computing and fuzzy approaches, e.g., \cite{Haque2026,Pervin2026,Ray2026}. Computing such steady-state density functions is in most cases only possible by extensive numerical simulations. One goal of this paper is to show how these calculations can be performed efficiently for a system with known existence of steady states and demonstrate the range of meaningful interpretations of the resulting probability densities.\medskip

Predator-prey systems are important for studying ecological and epidemiological dynamics \cite{GomezHernandez2024,Almatrafi2025a,Almatrafi2025b}. Specifically, the stochastic extensions of the Rosenzweig-MacArthur predator-prey ODE model are an active area of research \cite{Stollenwerk2022,Wyse2022,Wang2025}, in order to investigate the stochastic influences on well-known deterministic systems. Mostly, these stochastic extensions are in the direction of introducing randomness during the dynamical evolution of the system, leading to stochastic differential equations (SDE). In this concise study, we are interested in a different type of stochastic extension, investigating the influence of stochastic knowledge about the parameters of the model, leading to RDEs (see \cite{Jornet2023} for more details on the differences from SDEs), for obtaining information about the non-trivial steady state of the system. For this purpose, we utilize the computational approach presented recently \cite{Hoegele2026} and show how it can be adapted to determine steady-state densities.\medskip

The idea of \textit{quantum-like} or \textit{quantum-inspired modeling} is of growing interest in recent research \cite{Wichert2021,Khrennikov2023}, and this field is developing rapidly. In this paper, we adopt the general paradigm of quantum-like modeling in the sense that we apply some of the characteristics of quantum theory -- such as the superposition of parameter states within a system -- to non-quantum (e.g., macroscopic) objects of study, such as population dynamics. We emphasize that we do not utilize the full mathematical apparatus of quantum theory and quantum probabilities for this purpose, e.g., investigating spectral theory of self-adjoint operators on Hilbert spaces. Instead, we adopt the concept of incoherent quantum superposition (i.e., without phase information and potential interference patterns) by modeling it with classical mixture model probability density functions which are directly interpretable. In this sense, \textit{superposition} is interpreted in this work as the simultaneous possibility of different values for each continuous parameter, which we describe mathematically by mixture model probabilities with separated modes. For a broader discussion of this interpretation applied to the predator-prey system we refer to the discussion Section \ref{sec:discussion}. To our knowledge, this is a novel perspective on modeling superposition for parameters of a dynamical system in the broader meaning of quantum-like modeling (e.g., see \cite{Facco2022} for a discussion about the use of this term), allowing direct applications and conclusions utilizing classical probability theory.\medskip

\subsection{Motivation and related approaches}
This work is motivated by the integration of RDEs with parameter values represented by probabilistic mixture models and by the steady-state analysis of the resulting systems within a novel computational framework. This study focuses on both the direct calculation of steady-state distributions and their stability analysis. The latter is based on the eigenvalues of the Jacobian of the system, introducing a new quantity $\kappa$ for the probabilistic stability. Consequently, the general approach presented in Section \ref{sec:GeneralApproach} provides a framework for steady-state analysis of RDEs. \medskip

To our knowledge, the literature does not provide such a framework for parameter uncertainties and superpositions in RDEs. In particular, the computational strategy is based on a recent publication utilizing a likelihood/posterior formulation for random equations \cite{Hoegele2026}. The well-known nonlinear Rosenzweig-MacArthur predator-prey system serves as an illustrative example, which allows direct interpretability. Other current research close to this approach follows a similar line of argumentation since it contains the extension of a deterministic ODE to an RDE, focuses on the steady-state / equilibrium densities and performs illustrative simulation studies \cite{Bevia2024,Cortes2022}. The differences are that these studies do not provide a general framework based on the posterior densities of the solution, show no extension to mixture model parameters and are limited in their presentation to linear systems.

\section{Methods}

\subsection{General approach}
\label{sec:GeneralApproach}

\subsubsection{Classical presentation of steady states and stability of ODE systems}

The standard approach to compute steady states of an autonomous ODE system is first to reformulate the system as an explicit ODE system of first order, i.e., $\dot{\boldsymbol{x}} = \boldsymbol{f}(\boldsymbol{x})$ ($\boldsymbol{x} : \mathbb{R} \rightarrow \mathbb{R}^n$ and $\boldsymbol{f} : \mathbb{R}^n \rightarrow \mathbb{R}^n$) and second to solve $\boldsymbol{f}(\boldsymbol{x})=\boldsymbol{0}$. In general, solving such a (nonlinear) steady-state equation system is a difficult task, since its solvability is unclear (dynamical systems do not need to have steady states), and if solutions exist, how many there are and how they can be calculated (analytically or numerically) remain to be investigated.\medskip

Assuming that solutions $\boldsymbol{x}^{\star}$ exist of the steady-state equation system, i.e., $\boldsymbol{f}(\boldsymbol{x}^{\star})=\boldsymbol{0}$, the classical way to investigate their local stability is to determine the eigenvalues of the linearized version of the system, i.e., we calculate the Jacobian of the right-hand side $J_{\boldsymbol{f}}$ and evaluate the Jacobian at the steady state  $J_{\boldsymbol{f}}(\boldsymbol{x}^{\star})$. Local stability is determined by the eigenvalue distribution of this Jacobian, leading to the following standard argument: If all real parts of the eigenvalues of $J_{\boldsymbol{f}}(\boldsymbol{x}^{\star})\in\mathbb{R}^{n\times n}$ are negative, the steady state is \textit{asymptotically stable}, i.e., the solutions of the nonlinear ODE systems converge to such stable steady states. If the real part of at least one eigenvalue is positive, the solutions are \textit{unstable}, i.e., nearby solutions will diverge from that steady state after arbitrarily small perturbations \cite{Aulbach2004}. For the cases with real parts equal to zero, further analysis is necessary for nonlinear $\boldsymbol{f}$, which will be of no interest in this work. In order to get the eigenvalues $\lambda\in\mathbb{C}$ of a real matrix, we need to solve the characteristic polynomial $\text{det}(J_{\boldsymbol{f}}(\boldsymbol{x}^{\star}) - \lambda\,I) = 0$, with $I$ the identity matrix.\medskip

In summary, for the investigation of steady states and their stability, one needs to solve the steady-state equation system $\boldsymbol{f}(\boldsymbol{x})=\boldsymbol{0}$ for finding steady states $\boldsymbol{x}^{\star}$, and for the stability investigation one needs to solve $\text{det}(J_{\boldsymbol{f}}(\boldsymbol{x}^{\star}) - \lambda\,I) = 0$ in order to get the $n$ eigenvalues $\lambda_1,\dots,\lambda_n\in\mathbb{C}$ and consider the signs of their real parts.

\subsubsection{General extension to random variable parameters for RDE systems}
\label{sec:GenRandom}

In the following, we treat the right-hand side of the ODE system as being dependent on parameter random variables $A_1,\dots,A_k$, or short as depending on a random variable vector $\boldsymbol{A}$, i.e.,
\begin{align*}
\dot{\boldsymbol{x}} = \boldsymbol{f}(\boldsymbol{x};\boldsymbol{A})\;.
\end{align*}
Finding the steady states (depending now on the densities in $\boldsymbol{A}$) is equivalent to finding the probability distribution of the solutions of 
\begin{align*}
\boldsymbol{f}(\boldsymbol{x};\boldsymbol{A})=\boldsymbol{0}\;.
\end{align*}
This modification of the steady-state equations turns them into a (nonlinear) random equation system. A novel algorithm was recently presented to calculate the likelihood and posterior density of the solution space for general random equations \cite{Hoegele2026}. Applying this framework to the steady-state problem of RDEs, one needs to identify the general random equation system $\boldsymbol{M}(\boldsymbol{x};\boldsymbol{A}) = \boldsymbol{B} $ with independent random variable vectors $\boldsymbol{A}$ and $\boldsymbol{B}$ in our specific problem. We define $\boldsymbol{M}_{\text{steady}}(\boldsymbol{x};\boldsymbol{A}) :=\boldsymbol{f}(\boldsymbol{x};\boldsymbol{A})$, and let $\boldsymbol{B}$ be a random variable vector with zero mean and a small standard deviation compared to the other random variables. As a result, one gets an approximation of the posterior steady-state density function $\pi_{\text{steady}}(\boldsymbol{x})$ for the solution of this random equation system (utilizing non- or weakly informative priors).\medskip

For the stability analysis, we need to calculate the Jacobian of the right-hand side of the ODE system (which now also depends on the parameters $\boldsymbol{A}$), i.e., $J_{\boldsymbol{f}}(\boldsymbol{x};\boldsymbol{A})$ and investigate the density of the (in general complex) eigenvalues with the roots of the characteristic polynomial, i.e., $\text{det}(J_{\boldsymbol{f}}(\boldsymbol{x};\boldsymbol{A}) - \lambda\,I) = 0$. This is a polynomial in $\lambda$ and depends on the possible steady-state location $\boldsymbol{x}$ and the parameter densities $\boldsymbol{A}$. For simplicity, we transfer the search of the complex roots $\lambda\in\mathbb{C}$ to an equivalent search in $\mathbb{R}^2$ by identifying $\lambda = \sigma_1 + i\cdot \sigma_2$ with $\boldsymbol{\sigma}\in\mathbb{R}^2$. Then the search for the complex eigenvalues is identical to 
\begin{align}
M_{\text{eig},1}(\boldsymbol{x};\boldsymbol{A}) &:= \text{Re}\left\lbrace \text{det}(J_{\boldsymbol{f}}(\boldsymbol{x};\boldsymbol{A}) - (\sigma_1+i\cdot\sigma_2)\,I)\right\rbrace\label{Equ:pi_eig1}\\
M_{\text{eig},2}(\boldsymbol{x};\boldsymbol{A}) &:= \text{Im}\left\lbrace \text{det}(J_{\boldsymbol{f}}(\boldsymbol{x};\boldsymbol{A}) - (\sigma_1+i\cdot\sigma_2)\,I)\right\rbrace\;,\label{Equ:pi_eig2}
\end{align}
with the solution of the random equation system $\boldsymbol{M}_{\text{eig}}(\boldsymbol{x};\boldsymbol{A})=\boldsymbol{C}$ with $\boldsymbol{C}$ a two-dimensional random variable vector with zero mean and a small standard deviation compared to the other random variables (playing a role analogous to that of $\boldsymbol{B}$ for the steady-state density). This leads to the two-dimensional approximation of the posterior density function of the eigenvalues $\pi_{\text{eig}}(\boldsymbol{\sigma};\boldsymbol{x})$ depending on a given location $\boldsymbol{x}$.\medskip

For the stability analysis, we are interested in high-density regions of $\pi_{\text{steady}}(\boldsymbol{x})$, since they correspond to steady-state regions, and use $\pi_{\text{eig}}(\boldsymbol{\sigma};\boldsymbol{x})$ to compute the corresponding eigenvalue density. The stability is assessed by integrating the probability densities over negative real parts of the eigenvalues (i.e., $\sigma_1 <0$) and we define the new quantity
\begin{align}
\kappa(\boldsymbol{x}) := \int\limits_{-\infty}^{0} \int\limits_{-\infty}^{\infty } \pi_{\text{eig}}(\boldsymbol{\sigma};\boldsymbol{x})\;\text{d}\sigma_2\,\text{d}\sigma_1\;,\label{equ:kappa}
\end{align}
which summarizes the probability of all eigenvalues having a negative real part. Plotting the densities of $\pi_{\text{steady}}(\boldsymbol{x})$ and $\kappa(\boldsymbol{x})$ next to each other allows a straightforward interpretation: For each $\boldsymbol{x}$ with a high density for steady states in $\pi_{\text{steady}}(\boldsymbol{x})$ we determine the probabilistic stability by evaluating the same location in $\kappa(\boldsymbol{x})$ (only values very close to $1$ indicate asymptotic stability).\medskip

As we have seen, the main part of the modeling consists of identifying random equations of type $\boldsymbol{M}(\boldsymbol{x};\boldsymbol{A}) = \boldsymbol{B} $ in our steady-state analysis. As presented in \cite{Hoegele2026}, in order to numerically approximate the posterior probability density $\pi(\boldsymbol{x})$ of the solution space of such random equations (utilizing a non-informative prior), one needs to approximate an integral by Monte Carlo methods, specifically
\begin{align}
\pi(\boldsymbol{x})\propto \;  \frac{1}{N} \sum\limits_{n=1}^N\left(\prod\limits_{r=1}^R f_{B_r}(M_r(\boldsymbol{x};\boldsymbol{s}_{n}))\right)\;,\label{posterior}
\end{align}
where $R$ is the number of equations, $N$ is the number of samples, and $\boldsymbol{s}_n=(A_{1,n},\dots,A_{k,n})$ for $n=1,\dots,N$ are the sample vectors (i.e., the parameter sets independently sampled according to their densities). Since we assume that there is no interfering prior knowledge about the solution space, one needs only to normalize the numerically calculated right-hand side to the total probability mass of $1$ to get an approximation of the posterior density. 
Eventually, this results in the computational framework and algorithm for calculating $\pi_{\text{steady}}(\boldsymbol{x})$ as presented in Algorithm 1. For the calculation of the stability map $\kappa(\boldsymbol{x})$ a variation needs to be performed, which is presented in Algorithm 2.

\begin{table}[htbp]
\hrulefill\\
\textbf{Algorithm 1} (modeling and computation of $\pi_{\text{steady}}(\boldsymbol{x})$):
\begin{itemize}
\item[1.] Reformulate the dynamical ODE system into the form $\dot{\boldsymbol{x}} = \boldsymbol{f}(\boldsymbol{x};\boldsymbol{A})$, i.e., identify $\boldsymbol{f}$, select the parameter probability densities summarized in the random variable vector $\boldsymbol{A}$ and choose a $\boldsymbol{B}$, e.g., as a normal distribution with a small standard deviation compared to the densities in $\boldsymbol{A}$.
\item[2.] Draw $N$ samples $\boldsymbol{s}_n$ from $\boldsymbol{A}$.
\item[3.] Evaluate Equation (\ref{posterior}) with $\boldsymbol{M}_{\text{steady}}$, i.e., $R=n$, for each given phase space vector $\boldsymbol{x}$. Repeat this evaluation on a grid for $\boldsymbol{x}\in\mathbb{R}^n$ that includes all regions of non-negligible posterior density. 
\item[4.] Normalize this evaluated intensity map to probability mass $1$ in order to achieve an approximation of $\pi_{\text{steady}}(\boldsymbol{x})$.
\end{itemize}
\hrulefill
\end{table}

\begin{table}[htbp]
\hrulefill\\
\textbf{Algorithm 2} (modeling and computation of $\kappa(\boldsymbol{x})$):
\begin{itemize}
\item[1.] In addition to performing Step 1 of Algorithm 1, compute the Jacobian $J_{\boldsymbol{f}}$ and, based on it, the functions $M_{\text{eig},1}$ and $M_{\text{eig},2}$ given in Equations (\ref{Equ:pi_eig1}) and (\ref{Equ:pi_eig2}), and define $\boldsymbol{C}$ analogously to $\boldsymbol{B}$.
\item[2.] Draw $N$ samples $\boldsymbol{s}_n$ from $\boldsymbol{A}$.
\item[3.] Select a phase space vector $\boldsymbol{x}$ and evaluate Equation (\ref{posterior}) with $\boldsymbol{M}_{\text{eig}}$, i.e., $R=2$, for each complex number presented as real and imaginary part in the vector $\boldsymbol{\sigma}$. Repeat this evaluation on a grid for $\boldsymbol{\sigma}\in\mathbb{R}^2$ that includes all realistically possible eigenvalues. 
\item[4.] Normalize this evaluated intensity map to probability mass $1$ in order to achieve an approximation of $\pi_{\text{eig}}(\boldsymbol{\sigma};\boldsymbol{x})$.
\item[5.] Compute $\kappa(\boldsymbol{x})$ based on Equation (\ref{equ:kappa}) over the $\boldsymbol{\sigma}$ grid.
\item[6.] Repeat steps 3-5 for calculating $\kappa(\boldsymbol{x})$ in order to cover all regions in $\mathbb{R}^n$ with high probability density values of $\pi_{\text{steady}}(\boldsymbol{x})$.
\end{itemize}
\hrulefill
\end{table}

\subsection{The Rosenzweig-MacArthur model}

\subsubsection{The deterministic steady-state analysis}
\label{sec:PredPrey_Deterministic}

The normalized Rosenzweig-MacArthur predator-prey model with Holling type II functional response is given by \cite{Wang2025}
\begin{align*}
\frac{\text{d}x}{\text{d}s} &= x\,\left(1-\frac{x}{k}\right)-\frac{m\,x\,y}{1+x}\\ 
\frac{\text{d}y}{\text{d}s} &= -c\,y +\frac{m\,x\,y}{1+x}\;,
\end{align*}
where $x$ is essentially the normalized prey population, $y$ the normalized predator population, $k$ the carrying capacity of the prey population, $m$ the reproduction rate of the predators, and $c$ the death rate of the predator population.\medskip
  
In order to find the steady states, we set $\frac{\text{d}x}{\text{d}s}=\frac{\text{d}y}{\text{d}s}=0$, i.e.,
\begin{align*}
0 &= x\,\left(1-\frac{x}{k}\right)-\frac{m\,x\,y}{1+x}\\ 
0 &= -c\,y +\frac{m\,x\,y}{1+x}\;.
\end{align*}

For all fixed positive values of $k,m,c$ two steady-state solutions are trivial, e.g., $(x,y)=(0,0)$ (both populations are extinct) and $(x,y)=(k,0)$ (the prey is at maximum capacity and the predators are extinct). Interestingly, there is also a coexistence steady state with both population sizes greater than zero. This deterministic equilibrium solution (for $m,c,k>0$ with $m>c$ and $k>\frac{c}{m-c}$) is given by
\begin{align}
x=\frac{c}{m-c},\;y=\frac{1}{m}\,\left(1+\frac{c}{m-c}\right)\left(1-\frac{c}{k\,(m-c)}\right)\;.\label{equ:deterministic_steadystate}
\end{align}

There are many extensions of this model in the literature, most notably including disease and infection terms of the populations, leading to even more steady states \cite{GomezHernandez2024}. In this concise study, we want to focus mainly on the single non-trivial steady state described in Equation (\ref{equ:deterministic_steadystate}) and investigate its behavior under stochastic influences. 

In order to further investigate the stability, one needs the Jacobian of the right-hand side of the system
\begin{align*}
J(x,y) = \left(\begin{array}{cc}  1-\frac{2\,x}{k} -\frac{m\,y}{(1+x)^2} & -\frac{m\,x}{1+x} \\
\frac{m\,y}{(1+x)^2} & -c + \frac{m\,x}{1+x}
\end{array}\right)\;.
\end{align*}
For example, inserting the first trivial steady state $(x,y)=(0,0)$ leads to a diagonal Jacobian matrix and the two real eigenvalues $\lambda_1=1$ and $\lambda_2=-c$, i.e., one positive and one negative, leading to an \textit{unstable} steady state. In general, eigenvalues are determined by finding the roots of the characteristic polynomial $\text{det}(J(x,y)-\lambda\,I)=0$, leading to the quadratic formula
\begin{align*}
0 &= \lambda^2 - \left( 1-\frac{2\,x}{k} - \frac{m\,y}{(1+x)^2} -c + \frac{m\,x}{1+x}\right)\cdot\lambda\nonumber\\
&\qquad + \left(  1-\frac{2\,x}{k} - \frac{m\,y}{(1+x)^2} \right)\cdot\left( -c + \frac{m\,x}{1+x} \right) + \frac{m^2\,x\,y}{(1+x)^3}\;.
\end{align*}
One can easily substitute Equation (\ref{equ:deterministic_steadystate}) into the characteristic polynomial and calculate the roots analytically. For the parameters $k=1,m=1,c=\frac{1}{4}$ (identical to the modes in the simulation in Figure \ref{fig:Res:sim00}) this leads to a non-trivial steady state at $(\frac{1}{3},\frac{8}{9})$ with eigenvalues $\lambda_{1,2} = -\frac{1}{12} \pm i\cdot \frac{\sqrt{17}}{12}$, i.e., both have negative real part leading to an \textit{asymptotically stable} non-trivial steady state.

\subsubsection{The steady-state analysis of the RDE}

In the following, we are interested in the extension of the Rosenzweig-MacArthur problem regarding the parameters $k,m,c$ as independent random variables following probability density functions and investigate the behavior of the non-trivial steady-state solution. This step can be interpreted as introducing uncertainty about the true values of these parameters. \medskip

Further, in the spirit of recent developments \cite{Hoegele2026}, we want to utilize multimodal mixture model densities for the parameters $m$ and $c$. This can be interpreted as the predator population being heterogeneous, meaning there are subpopulations for which different parameter values essentially apply, and we allow all combinations of these parameter values. The new perspective we want to propose in this work is that the mixture model represents the simultaneous presence of different values for the parameters, i.e., a superposition. A more detailed discussion of these interpretations (and especially the interpretation of the simulated results) is provided in the discussion in Section \ref{sec:discussion}.\medskip

In order to calculate the steady states, this is given by (see Section \ref{sec:GenRandom})
\begin{align*}
M_{\text{steady},1}(\boldsymbol{x};\boldsymbol{A}) & := x_1\,\left(1-\frac{x_1}{A_1}\right)-\frac{A_2\,x_1\,x_2}{1+x_1}\\
M_{\text{steady},2}(\boldsymbol{x};\boldsymbol{A}) & := -A_3\,x_2 +\frac{A_2\,x_1\,x_2}{1+x_1}
\end{align*}
and $x_1:= x$, $x_2:=y$ as well as the random variables $A_1:=k\sim f_k$, $A_2:=m\sim f_m$ and $A_3:=c\sim f_c$. \medskip

In order to analyze the stability of the steady states, we will investigate the eigenvalues of the Jacobian. As demonstrated before, we search for solutions $\lambda\in\mathbb{C}$ of (with the same $A_1,A_2,A_3$ as before)
\begin{align*}
0 &= \lambda^2 - \left( 1-\frac{2\,x_1}{A_1} - \frac{A_2\,x_2}{(1+x_1)^2} -A_3 + \frac{A_2\,x_1}{1+x_1}\right)\cdot\lambda\nonumber\\
&\qquad + \left(  1-\frac{2\,x_1}{A_1} - \frac{A_2\,x_2}{(1+x_1)^2} \right)\cdot\left( -A_3 + \frac{A_2\,x_1}{1+x_1} \right) + \frac{A_2^2\,x_1\,x_2}{(1+x_1)^3}\;.
\end{align*}
Since $\lambda$ is in general complex, we identify it with its real and imaginary parts, i.e., two real values $\lambda = \sigma_1 + i\cdot \sigma_2$ with $\boldsymbol{\sigma}\in\mathbb{R}^2$, which leads to (see Section \ref{sec:GenRandom})
\begin{align*}
0 &= \sigma_1^2-\sigma_2^2 + i\cdot 2\sigma_1\sigma_2 - \left( 1-\frac{2\,x_1}{A_1} -  \frac{A_2\,x_2}{(1+x_1)^2} -A_3 + \frac{A_2\,x_1}{1+x_1}\right)\cdot(\sigma_1 + i\cdot \sigma_2)\nonumber\\ 
&\qquad + \left(  1-\frac{2\,x_1}{A_1} -  \frac{A_2\,x_2}{(1+x_1)^2} \right)\cdot\left( -A_3 + \frac{A_2\,x_1}{1+x_1} \right) + \frac{A_2^2\,x_1\,x_2}{(1+x_1)^3}\;,
\end{align*}
which will be further interpreted as forcing the real and imaginary part of the right-hand side equal to zero, i.e.,
\begin{align*}
M_{\text{eig},1}(\boldsymbol{\sigma};\boldsymbol{x},\boldsymbol{A}) & := \sigma_1^2-\sigma_2^2 - \left( 1-\frac{2\,x_1}{A_1} -  \frac{A_2\,x_2}{(1+x_1)^2} -A_3 + \frac{A_2\,x_1}{1+x_1}\right)\cdot\sigma_1\nonumber \\
&\qquad + \left(  1-\frac{2\,x_1}{A_1} - \frac{A_2\,x_2}{(1+x_1)^2} \right)\cdot\left( -A_3 + \frac{A_2\,x_1}{1+x_1} \right) + \frac{A_2^2\,x_1\,x_2}{(1+x_1)^3} \\
M_{\text{eig},2}(\boldsymbol{\sigma};\boldsymbol{x},\boldsymbol{A}) &:= 2\sigma_1\sigma_2 - \left( 1-\frac{2\,x_1}{A_1} -  \frac{A_2\,x_2}{(1+x_1)^2} -A_3 + \frac{A_2\,x_1}{1+x_1}\right)\cdot \sigma_2\;.
\end{align*}
The corresponding posterior density function of the eigenvalues $\pi_{\text{eig}}(\boldsymbol{\sigma};\boldsymbol{x})$ depending on a given $\boldsymbol{x}$ can again be calculated by Monte Carlo simulation by Equation (\ref{posterior}). Based on this, the stability probability $\kappa(\boldsymbol{x})$ can be calculated according to Equation (\ref{equ:kappa}).\medskip

\subsection{Independent verification}
\label{sec:IndepVerif}

There are three classical ways to verify the steady-state computation. All three include random sampling parameter sets from the described distributions and use them in the deterministic system in order to approximate the steady-state density. In the spirit of Monte Carlo simulations, one then needs to repeat this until a 2D histogram of steady states is obtained. The three approaches are: 
\begin{itemize}
\item[a)] Using the analytically calculated equilibrium solution if this is possible. This would be the numerically most efficient way, but we want to clarify that such analytic solutions are rarely available for arbitrary ODEs. Due to this, we deliberately avoid using this verification approach in this study.
\item[b)] Using the deterministic nonlinear steady-state equation system. This leads to essentially solving a random equation utilizing repeatedly an iterative nonlinear equation solver. This approach does not need analytical solvability and can be formulated generally, but it is numerically more expensive than a) and relies directly on the solver algorithm.
\item[c)] Using the ODE as the deterministic system. This approach repeatedly uses an ODE solver for each parameter set to approximate the RDE-derived steady-state distribution. This can also be formulated very generally, but it is numerically the most expensive approach and now depends on the ODE solver algorithm and the initial values of the system. Specifically, if one does not start the simulation directly at an unstable steady state, this approach will naturally find only stable steady states.
\end{itemize}

\section{Results}

\subsection{Steady-state analysis}
In the results section, we examine several simulation results for the non-trivial steady-state distributions based on the described random equations. The following plots show on the left the parameter probability densities presented as sample histograms utilized in the Monte Carlo integrations of the calculation for $k$, $m$ and $c$ (from top to bottom) (see \cite{Hoegele2026} for details). We utilize Gaussian mixture models for simplicity without restricting the applicability of the general methodology. On the right side of the plots, the resulting steady-state posterior distribution is presented showing the color-coded probability of the steady state for $x_1$ and $x_2$ as coordinates. The number of samples in the Monte Carlo integration is $N=24000$ (see Equation (\ref{posterior})) generated by Latin hypercube sampling. This sample number was chosen since the computational cost is comparably low, the histograms approximate all density functions well (as demonstrated in the results figures on the left side) and the resulting posterior densities $\pi_{\text{steady}}$ show no relevant numerical noise. The standard deviations of $f_{B_1}$ and $f_{B_2}$ are with $\sigma_{B_1}=\sigma_{B_2}=0.005$ significantly smaller than all other standard deviations in the simulation (which can directly be compared in the figures), and it is verified in Figure \ref{fig:Res:sim03_verif} that this is small enough to not compromise the obtained results.\medskip

In general, a compromise for numerical efficiency has to be found -- taking standard deviations for $\boldsymbol{B}$ that are too small leads in Equation (\ref{posterior}) to more evaluations of the density of $\boldsymbol{B}$ at locations $M_r(\boldsymbol{x};\boldsymbol{s}_{n})$ where the density is zero, i.e., it does not contribute to the final density function but generates calculation costs, and the posterior estimate may become noisy (eventually leading to the need for a higher number of samples in Monte Carlo integration). If the standard deviation is too large, the result may be substantially compromised, leading to blurred versions of the posterior density. A good compromise can be found by starting with relatively large standard deviations and reducing it step by step (by possibly simultaneously increasing sample sizes due to noise) until no effective difference in the final posterior is visible. This has been done for the following simulations.  \medskip

\begin{figure}[htbp]
\centering
\includegraphics[width=16cm]{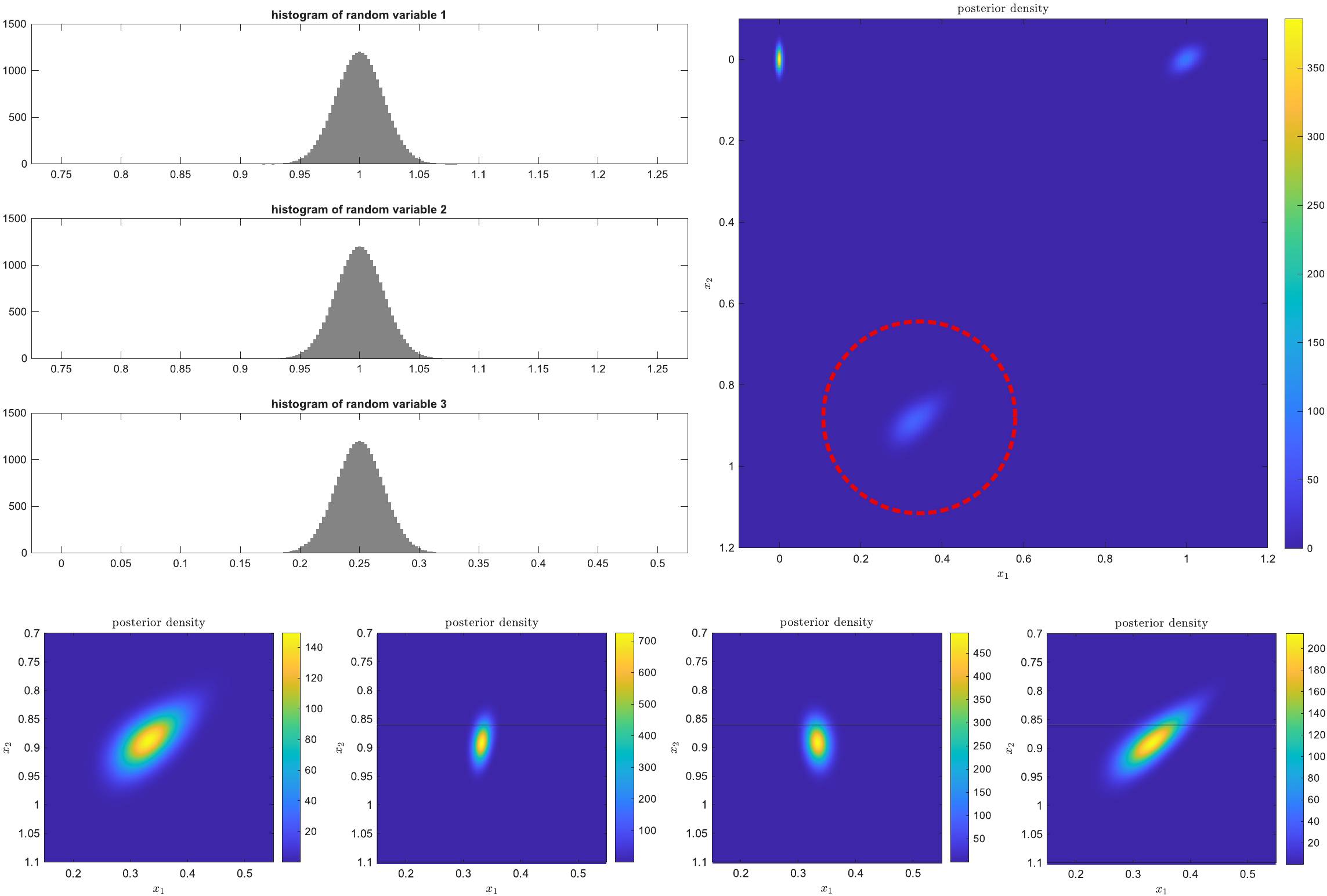}
\caption{Top row: Left: Unimodal parameter densities containing broad uncertainties about the parameters. Right: $\pi_{\text{steady}}(\boldsymbol{x})$ as intensity plot. The three expected steady states (the top two are the trivial solutions) and the red-circled density region corresponds to the non-trivial equilibrium state. The observed shape of this density shows how the parameter uncertainties propagate to the equilibrium steady state and is non-trivial due to the nonlinearity of the equations. Bottom row (from left to right): Sensitivity analysis of the parameters on the non-trivial steady state. First plot: Zoomed version of the non-trivial steady state with the unimodal parameter densities as in the top row. Second plot: only $k$ is set to the presented density, $m$ and $c$ fixed values at the modes. Third plot: only $m$ is set to the presented density, $k$ and $c$ fixed values at the modes. Fourth plot: only $c$ is set to the presented density, $k$ and $m$ fixed values at the modes. }  
\label{fig:Res:sim00}
\end{figure}

\begin{figure}[htbp]
\centering
\includegraphics[width=16cm]{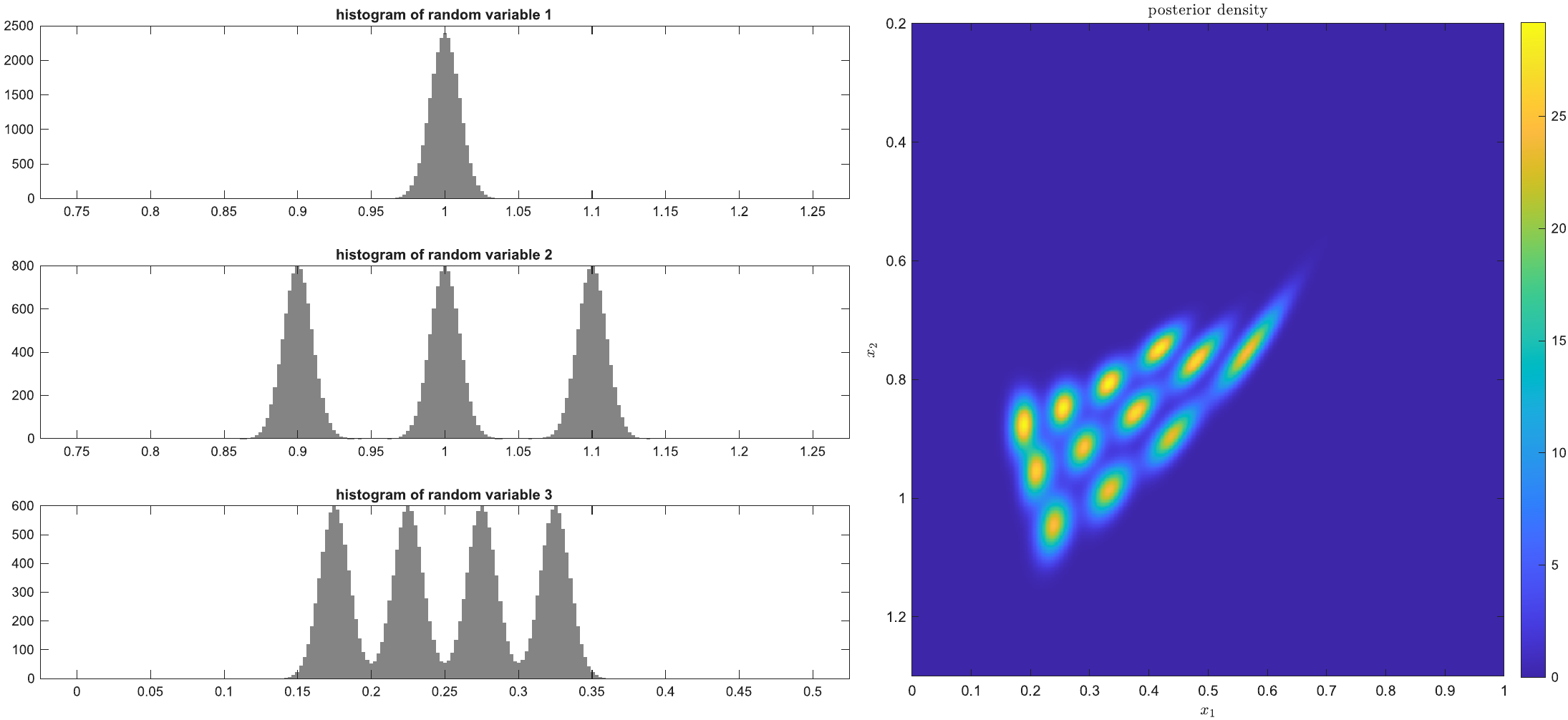}
\caption{Left: The parameters $m$ and $c$ (second and third row) show distinct and even separated multimodal patterns. Right:  $\pi_{\text{steady}}(\boldsymbol{x})$ as intensity plot. This leads to the non-trivial steady-state distribution with $12$ peaks, which comes from $3$ modes for $m$ and $4$ modes for $c$ ($3\times 4=12$ parameter combinations). }  
\label{fig:Res:sim01}
\end{figure}

\begin{figure}[htbp]
\centering
\includegraphics[width=16cm]{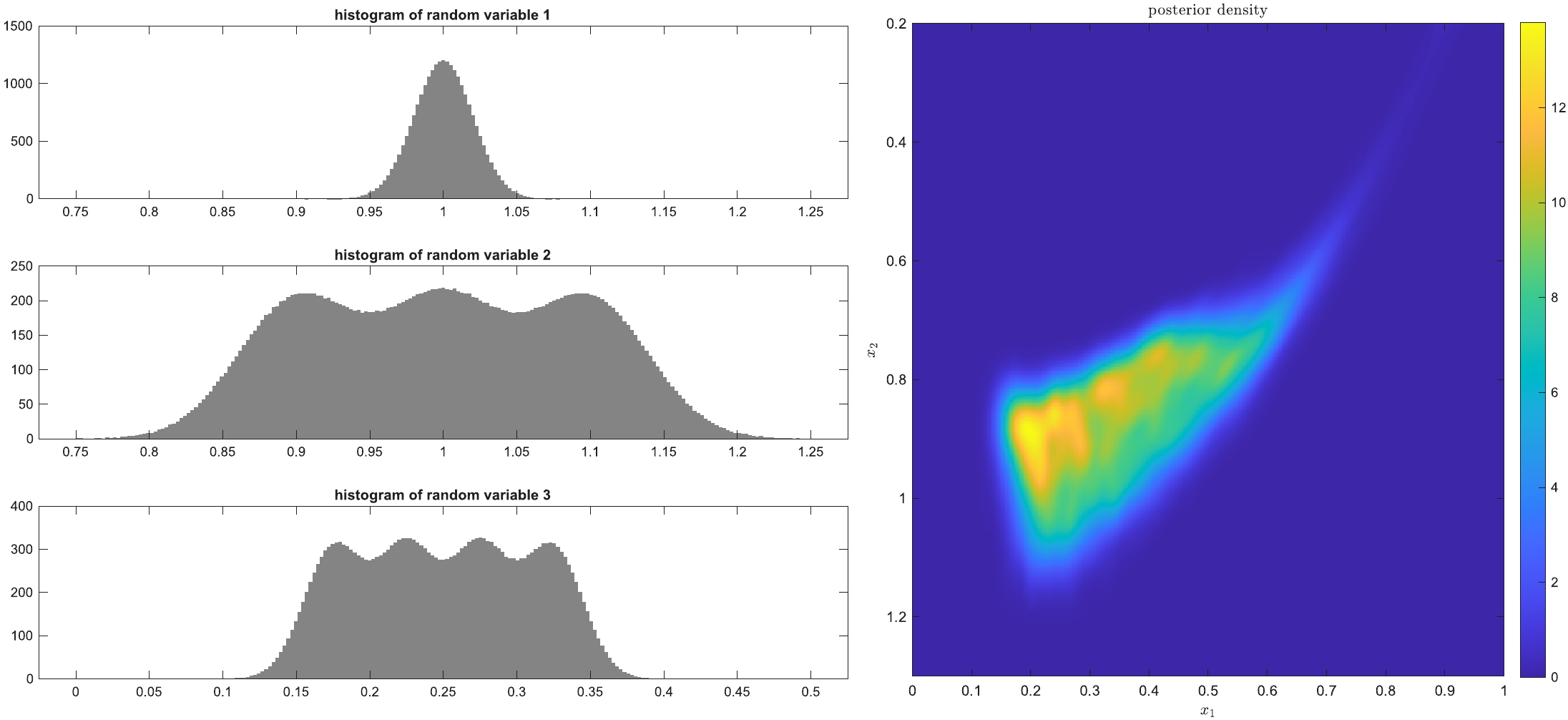}
\caption{Left: The parameters $m$ and $c$ are again multimodal, but without clearly separated modes, resulting in broad non-distinct parameter densities. Right:  $\pi_{\text{steady}}(\boldsymbol{x})$ as intensity plot. The non-trivial steady-state distribution shows a large blurred region of characteristic shape where individual modes are hardly recognizable. }  
\label{fig:Res:sim02}
\end{figure}

\begin{figure}[htbp]
\centering
\includegraphics[width=16cm]{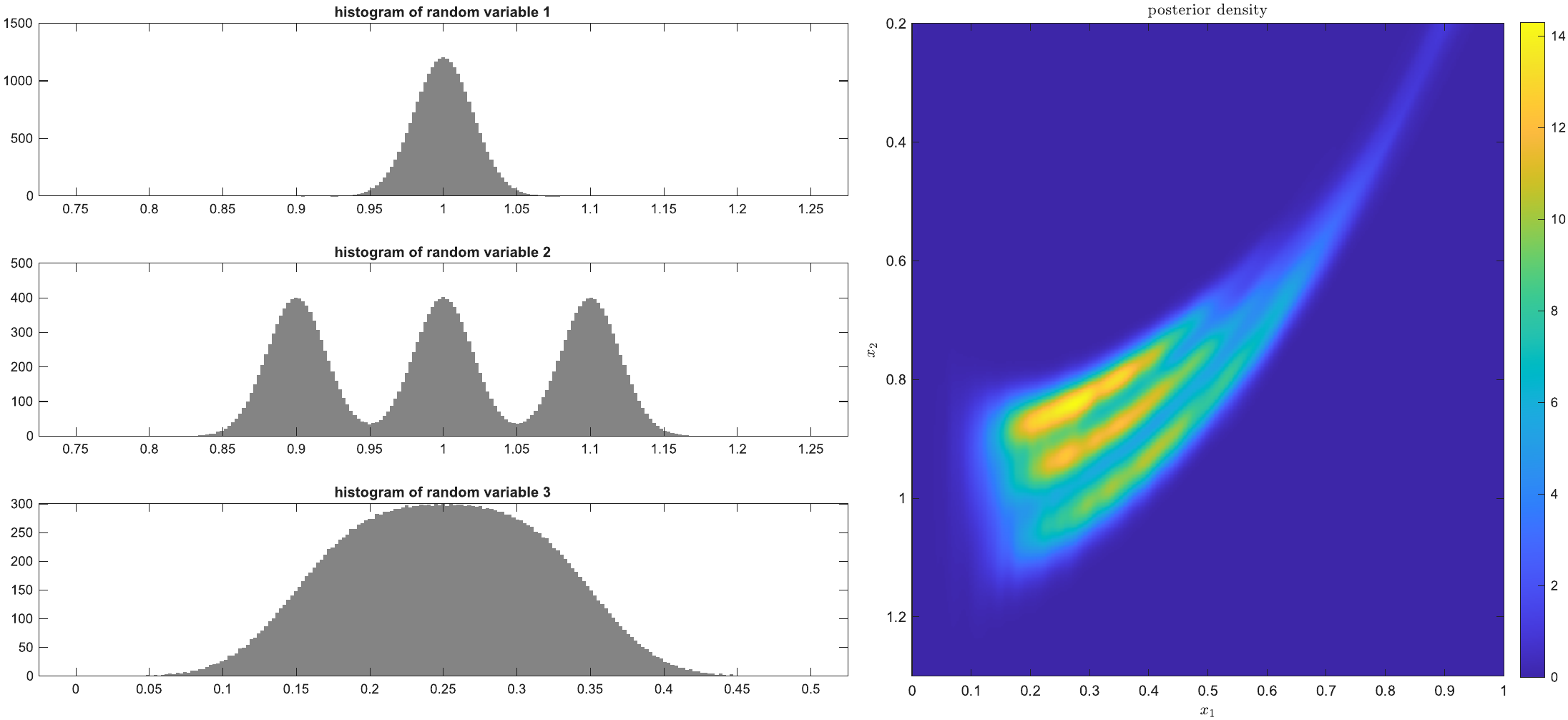}
\caption{Left: The parameter $m$ shows distinct modes and $c$ only a broad non-distinct probability distribution which is composed of four Gaussians. Right:  $\pi_{\text{steady}}(\boldsymbol{x})$ as intensity plot. The non-trivial steady-state distribution shows how this mixture of parameters propagates to the steady states, again showing non-trivial patterns due to nonlinearity.}  
\label{fig:Res:sim03}
\end{figure}

\begin{figure}[htbp]
\centering
\includegraphics[width=14cm]{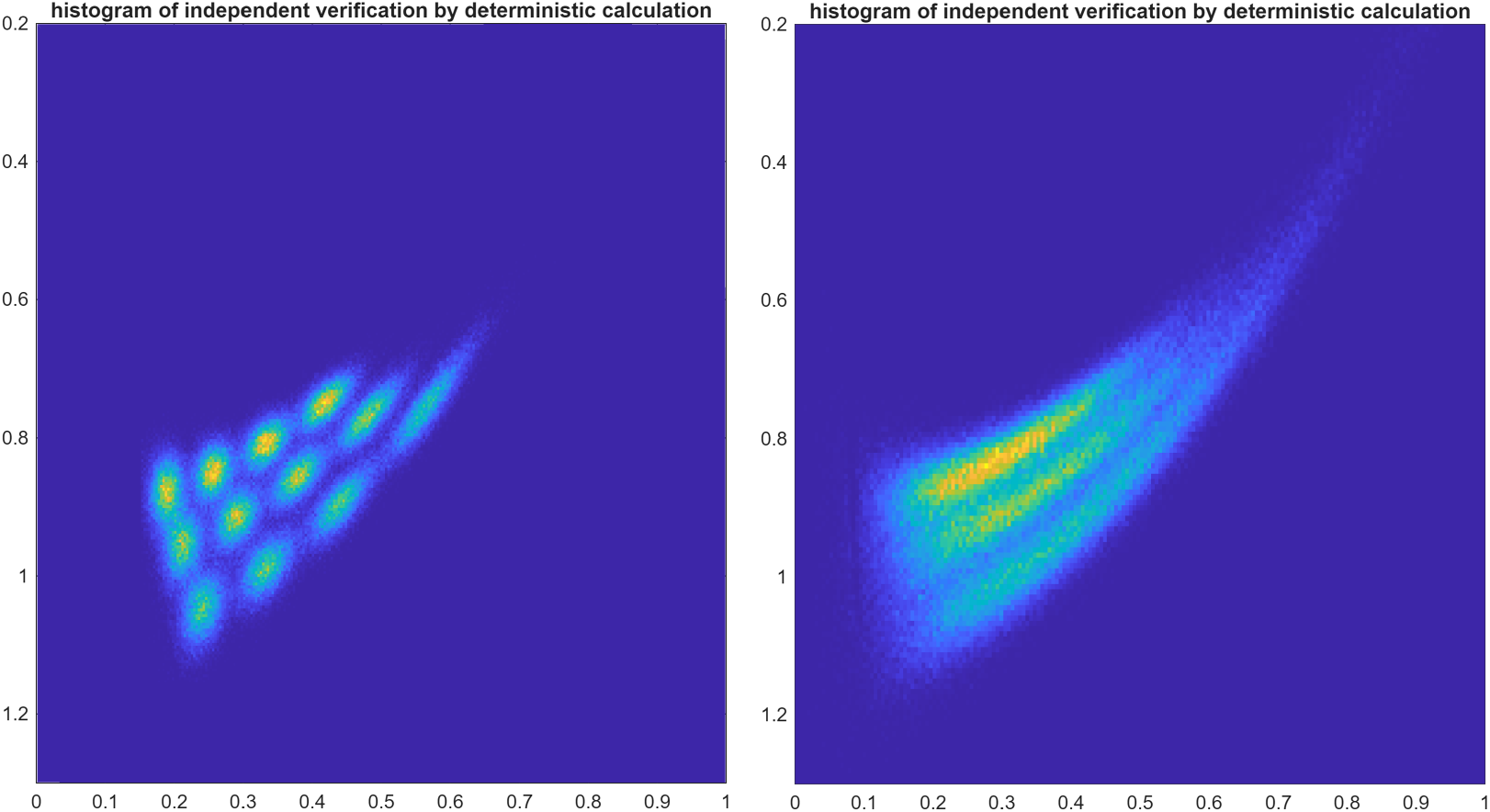}
\caption{2D histograms of the verification calculation of the simulation examples presented in Figure \ref{fig:Res:sim01} (left) and in Figure \ref{fig:Res:sim03} (right) utilizing $N=240000$ numerical solutions of the sampled deterministic steady-state equations as described in Section \ref{sec:IndepVerif}. MATLAB's fsolve iterative solver was utilized with the constant start value $x=[0.3,0.9]$.}  
\label{fig:Res:sim03_verif}
\end{figure}

Specifically in Figure \ref{fig:Res:sim00} (top row), we present the three steady-state solutions (two trivial, one non-trivial) for unimodal parameter densities. In all further plots, we only focus on the non-trivial steady state. Essentially, we get blurred versions of the deterministic solutions in Section \ref{sec:PredPrey_Deterministic} due to the introduced uncertainties. In Figure \ref{fig:Res:sim00} (bottom row) an additional sensitivity analysis on the non-trivial steady state is performed, showing the increasing sensitivity in the order of $k$, $m$ and $c$. It is additionally demonstrated that the introduced blur direction of the non-trivial steady-state distribution is different for each individual parameter. In Figure \ref{fig:Res:sim01}, we show the effect of distinctly separated modes in the mixture model for $m$ and $c$. This leads to a complex structure for the steady-state density, essentially containing all possible combinations of the modes. In Figure \ref{fig:Res:sim02}, we show the resulting steady-state distribution for non-distinct multimodal densities for $m$ and $c$. In Figure \ref{fig:Res:sim03}, the influence is presented when $m$ contains distinct and separated modes and $c$ does not. Overall, the interplay of the combinatorial complexity of mixture models in the parameters with the nonlinearity of the predator-prey equations leads to an emergence of complex structures for the non-trivial steady-state density.\medskip

For verification, in Figure \ref{fig:Res:sim03_verif}, the 2D histogram approximating the posterior is calculated independently by the verification approach presented in Section \ref{sec:IndepVerif} by solving the nonlinear steady-state equations numerically with an iterative solver for $N=240000$ drawn deterministic equations (approach (b)). It can be seen that (besides minor numerical artifacts) the posteriors are approximated well in their overall appearance, verifying the calculations based on Equation (\ref{posterior}). This shows that, although the simulations for Figure \ref{fig:Res:sim03_verif} take significantly longer compared to the results based on Equation (\ref{posterior}), the detailed structure of the steady-state density is not yet fully visible and would require even longer computation runtime.

\subsection{Stability analysis}

\begin{figure}[htbp]
\centering
\includegraphics[width=9cm]{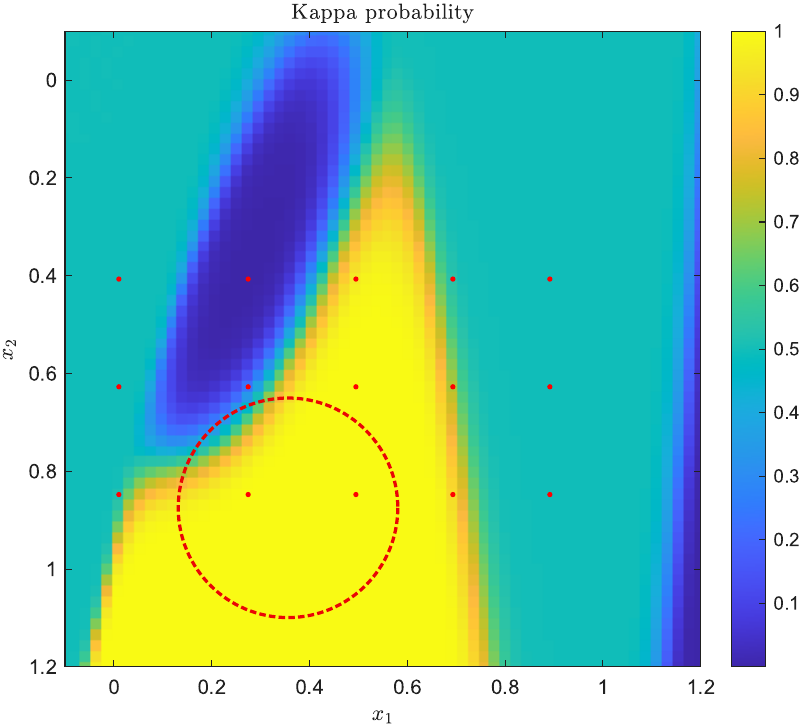}
\caption{Probabilistic stability map $\kappa(\boldsymbol{x})$ for the case presented in Figure \ref{fig:Res:sim00}. Values close to $1$ indicate a high probability of asymptotic stability. Further, a 15-point grid is presented for a selected case study of the corresponding eigenvalues. Further, the same red dashed circle as presented in Figure \ref{fig:Res:sim00} is presented in order to show the high-intensity region of the steady-state distribution.}  
\label{fig:Res:sim00_stability_kappa}
\end{figure}

\begin{figure}[htbp]
\centering
\includegraphics[width=16cm]{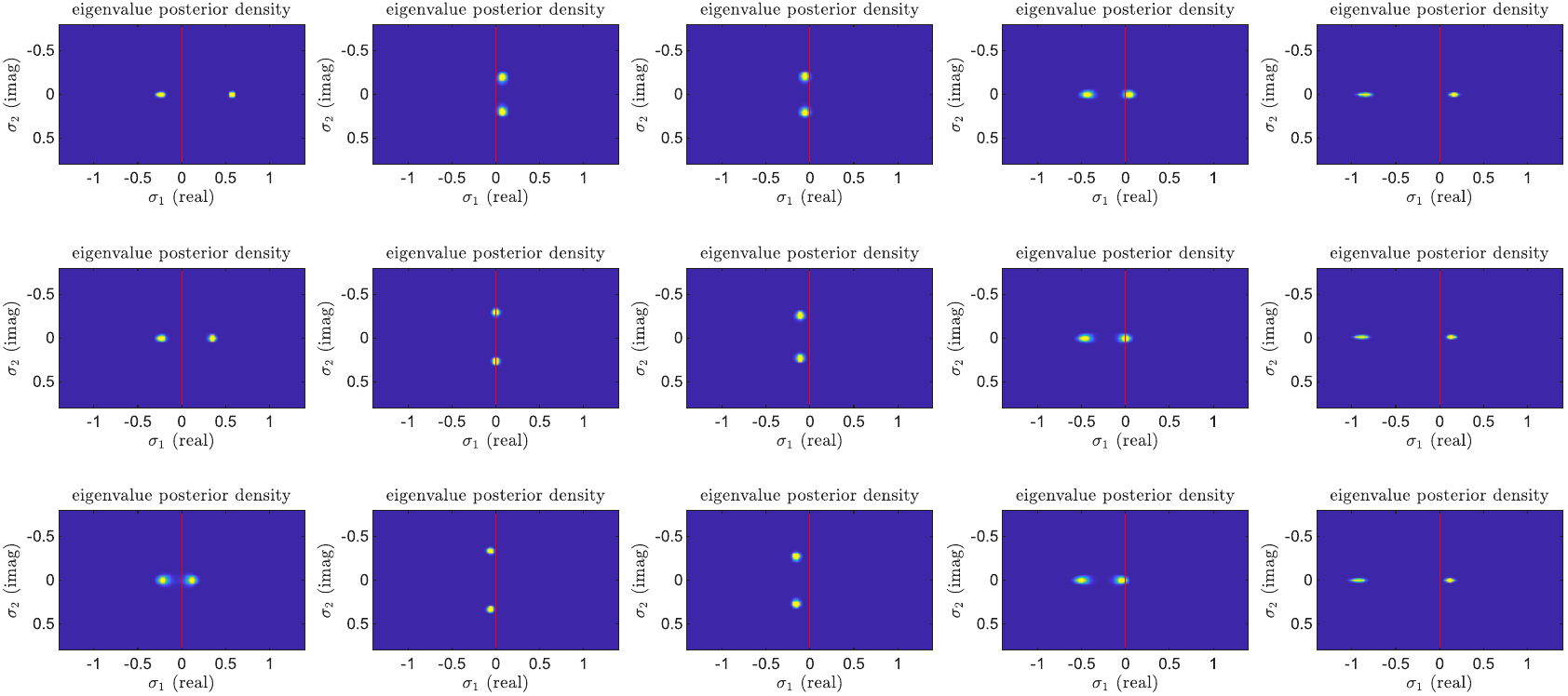}
\caption{Intensity plot of the eigenvalue densities $\pi_{\text{eig}}(\boldsymbol{\sigma};\boldsymbol{x})$ corresponding to the 15 point grid for $\boldsymbol{x}$ as presented in Figure \ref{fig:Res:sim00_stability_kappa}. If all intensities are on the left of the red line (negative real part of the eigenvalues, i.e., $\sigma_1 < 0$) asymptotic stability is guaranteed.}  
\label{fig:Res:sim00_stability_eigval}
\end{figure}

We perform the stability analysis for the first case, presented in Figure \ref{fig:Res:sim00}, and the second case, presented in Figure \ref{fig:Res:sim01}, as two distinct examples. In Figure \ref{fig:Res:sim00_stability_kappa} the probabilistic stability map is presented for the first case. It shows asymptotically stable regions with probability 1 (yellow) and unstable regions (all others) to different degrees. Since stability is only important in regions with high steady-state probabilities, the same dashed circle as shown in the steady-state density in Figure \ref{fig:Res:sim00} around the non-trivial steady state is presented. This shows that throughout the region of the non-trivial steady state asymptotic stability is guaranteed. Further, the other two trivial steady states around $(0,0)$ and $(k,0)$ are unstable, as expected. Every pixel in the $\kappa$ probability plot is calculated by evaluating the full eigenvalue posterior density. For a 15-point grid (presented as red dots), the corresponding eigenvalue densities are shown as intensity plots in Figure \ref{fig:Res:sim00_stability_eigval}. Since in this first case the probability densities of $A_i$ are all unimodal, we essentially expect to see blurry versions of the distinct eigenvalue positions. Note that, in this problem the characteristic polynomial contains only real-valued coefficients and is of degree 2. This means, by the fundamental theorem of algebra, that we have either only two real eigenvalues (on the real line, i.e., $\sigma_2=0$) or two complex conjugate eigenvalues (symmetric with respect to the $\sigma_1$-axis), which can be directly seen and which also serves as the verification of the correct calculation of the eigenvalue densities.

\begin{figure}[htbp]
\centering
\includegraphics[width=9cm]{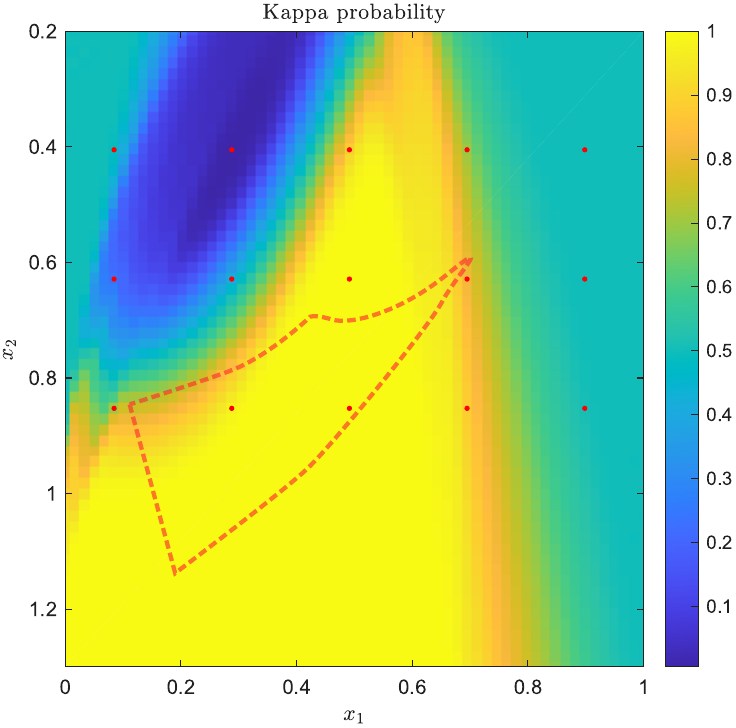}
\caption{Probabilistic stability map $\kappa(\boldsymbol{x})$ for the case presented in Figure \ref{fig:Res:sim01}. Values close to $1$ indicate a high probability of asymptotic stability. Further, the same 15 point grid as in Figure \ref{fig:Res:sim00_stability_kappa} is presented for a selected case study of the corresponding eigenvalues. Further, the high-intensity region of the steady-state distribution in Figure \ref{fig:Res:sim01} is presented as red dashed shape.}  
\label{fig:Res:sim01_stability_kappa}
\end{figure}

\begin{figure}[htbp]
\centering
\includegraphics[width=16cm]{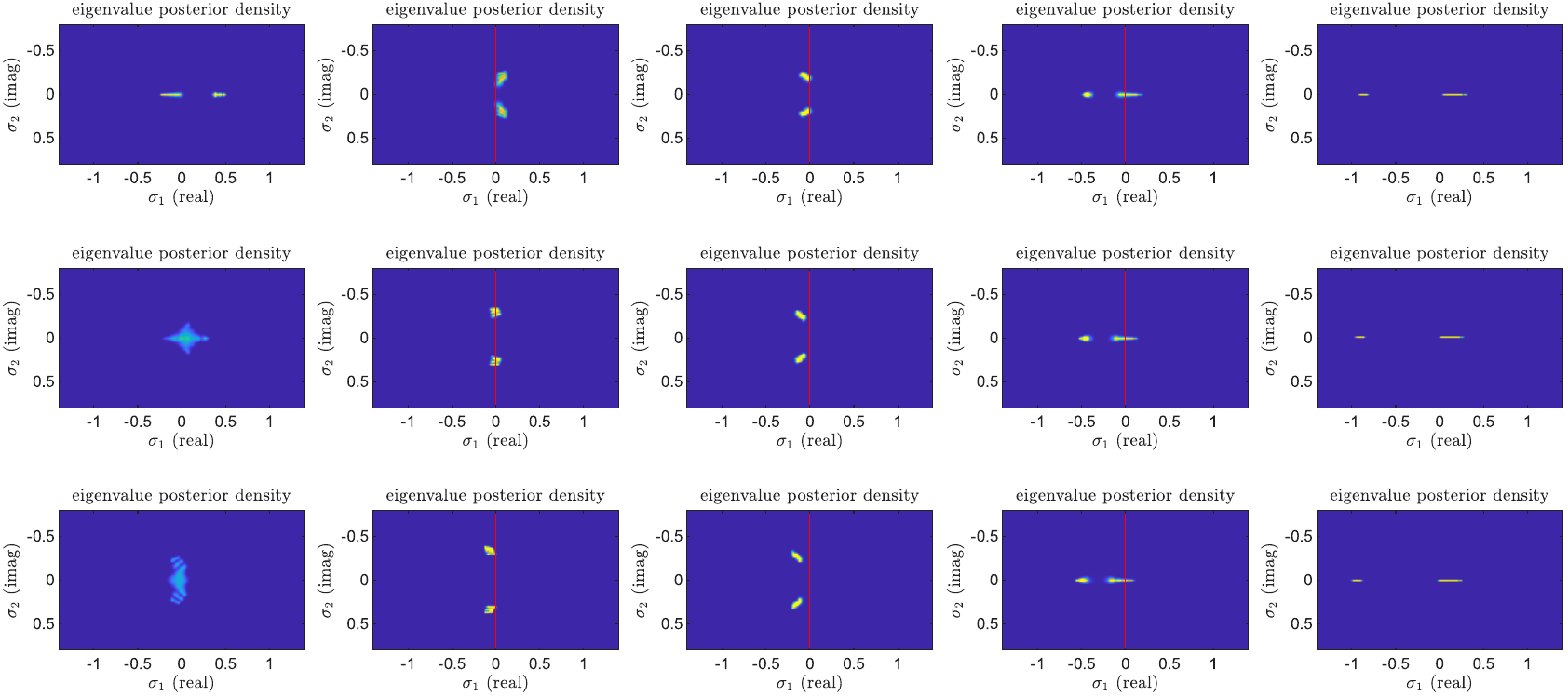}
\caption{Intensity plot of the eigenvalue densities $\pi_{\text{eig}}(\boldsymbol{\sigma};\boldsymbol{x})$ corresponding to the 15 point grid as presented in Figure \ref{fig:Res:sim01_stability_kappa}. If all intensities are on the left of the red line (negative real part of the eigenvalues, i.e., $\sigma_1 < 0$) asymptotic stability is guaranteed.}  
\label{fig:Res:sim01_stability_eigval}
\end{figure}

In Figure \ref{fig:Res:sim01_stability_kappa}, the $\kappa$ probability plot is presented for the second case, analogous to Figure \ref{fig:Res:sim00_stability_kappa} containing superpositions of the parameters. This time, the high-intensity steady-state region of Figure \ref{fig:Res:sim01} has a characteristic outer shape and is presented by the dashed line figure. It can be seen that the steady states are asymptotically stable throughout nearly the entire region. This means that, although the position of the non-trivial steady state is rather sensitive with respect to the parameter densities, the stability of these steady states remains robust. Further, the same grid for the eigenvalue evaluation as in Figure \ref{fig:Res:sim00_stability_kappa} is presented as red dots. In Figure \ref{fig:Res:sim01_stability_eigval}, we see the corresponding eigenvalue distribution, which shows more complex structures, since the $A_i$ have multimodal densities which generate also shaped eigenvalue densities. Still, the symmetry with respect to the $\sigma_1$-axis works as a verification of the numerical calculation.  

\section{Discussion}
\label{sec:discussion}
There are several interpretations of the steady-state density results for the predator-prey model listed below, which should be distinguished. Besides the demonstration of the numerical results, these interpretations are a major point of the study and should broaden its perspective.
\begin{itemize}
\item In the first interpretation, we do not know the specific parameter values; nevertheless, there is only one true parameter set active in the population. The fact that there are multimodal mixture models for the parameter densities then reflects the fact that one has different (possibly non-compatible) sources of information about these parameters. Then, the results show the probability distribution about the non-trivial steady state as a probability density function of the true parameter set. \textit{Remark:} Although there may be non-compatible sources of separated stochastic information about the possible parameter sets, the effective steady-state distribution can be overlapping, i.e., it may not necessarily show these separated possible combinations (e.g., as visible in Figure \ref{fig:Res:sim02}). 
\item Second, similar to the first interpretation, at each time and location of an investigated population, there is only one true parameter set, but either this parameter set varies in time and/or space, with potentially many neighboring (essentially non-interacting) subpopulations containing all combinations of parameter density modes. The calculated distribution then shows the overall (averaged) distribution of steady states related to the sampling of the populations for different times and locations. 
\item Third, we say that the populations are heterogeneous, and individual member behavior is determined by parameter values drawn from the parameter densities (not necessarily separated populations in time and/or space) and, therefore, are best represented by a mixture distribution over all scales of spatiotemporal averaging of the populations (without a hidden unknown true parameter set of whole subpopulations). Then the results show the expected steady-state population distribution, and for each complete draw of population parameters, we just observe one possible realization of that steady-state distribution. One could then argue that a characteristic value (such as the expectation value or the mode(s)) of the posterior density should be expected in observations on average. But this must be interpreted carefully, since it could be that the expectation value lies in a region where the posterior density is zero, e.g., for multimodal posteriors as presented in the examples.
\item Fourth, taking the third interpretation one step further, we assume that each member of the populations follows actions effectively according to the parameter set densities, i.e., we completely lose the assumption that individuals would stick to a single true parameter set while we as observers simply do not know which is the true one. This means, we assume the probability distributions as the only basic truth, so that each member is active simultaneously in different parameter states -- a superposition of the behavior. In this context, the resulting steady-state distribution as a whole represents the true steady state.   
\end{itemize}

The last interpretation is unconventional and is of particular interest, especially when complicated distribution patterns due to the nonlinearity of the problem occur, such as in the presented example. This is related to interpretations of quantum probability that have been applied to ecological systems where the basic truth is the distribution without a hidden true state of parameters \cite{Angeler2024, Huebsch2024, Bernardini2023}. Interestingly, one can then interpret the solution of the steady-state equations directly as a nonlinear operation acting on these variables in superposition. In this context, we are not concerned with the \textit{measurement} or \textit{collapse of the wave function} side of quantum theory, which could be interpreted accordingly as the counting of the populations at a specific point in time.\medskip 

Most importantly, these interpretations should not necessarily be considered as mutually exclusive, since they can be valid simultaneously or regarded as partially  complementary perspectives. For example, parts of the parameter probability distributions might come from uncertainty or missing knowledge, and other parts could come from inhomogeneity between individuals or true superposition of parameter values for individuals in the population. An important point is that these are only interpretations of the probabilities informing the assumed parameter densities. The mathematical treatment once these densities are established is straightforward by applying arguments of classical probability theory.\medskip

With respect to the stability analysis, it is demonstrated how the framework in \cite{Hoegele2026} can be utilized also to support a more detailed analysis of these steady-state densities. One major obstacle is that the solution of the characteristic polynomial leads to possibly complex eigenvalues, which is solved by setting simultaneously the real and imaginary parts of the characteristic polynomial to zero, leading to two real-valued random equations. Since this is general, both the steady-state intensity calculation and its stability analysis can be easily extended to arbitrary ODE systems.\medskip

From a mathematical modeling and numerical perspective, applying the probabilistic method in \cite{Hoegele2026} allows an efficient calculation of the steady states of random differential equations without the computationally expensive time-dependent simulation until a steady state is reached or the iterative solution of the nonlinear steady-state equation system as presented in Figure \ref{fig:Res:sim03_verif}. It is one goal of this study to demonstrate this efficiency and elegance for the presented example.\medskip

It is well known that findings for predator-prey models are often directly applicable to epidemiological dynamics \cite{GomezHernandez2024,Berkal2026} and in this context the calculated distributions can help to estimate public health needs when parts of the disease and/or recovery parameters are uncertain or show multimodality due to diversity in the population.\medskip

Further, it is evident that this framework to approximate steady-state distributions and their stability analysis can be directly applied to any ODE system that one wants to extend to a RDE, since all explicit autonomous ODE systems can be reformulated to first order, and by setting all time derivatives to zero the steady-state solutions can be obtained by solving the residual algebraic (in general nonlinear) random equation system as presented. One needs only to identify the elements of the random equation system $\boldsymbol{M}(\boldsymbol{x};\boldsymbol{A}) = \boldsymbol{B} $ and model the parameter densities. Of course, this also applies to higher-dimensional systems, such as predator-prey simulations with many different populations or other dynamical systems. For this purpose, future work should systematically investigate computationally efficient methods for evaluating the posterior in high-dimensional systems, e.g., with importance sampling or avoiding calculations in low-density regions. Thus, this study can also be read as an example of how to calculate steady states in other domains that combine dynamical behavior with uncertainties or superpositions. Another angle of future research could be how this approach can be extended to solving steady-state solutions for partial differential equations (PDE) with parameter uncertainties. For example, a direct connection occurs if the PDE is reformulated as a large ODE system in the course of numerical approximation.\medskip

\section{Conclusion}

In this concise study, we have shown how stochastic information about the parameters of the Rosenzweig-MacArthur predator-prey model leads to probability distributions of the non-trivial steady state. Further, we show that this extension leads to complex steady-state distributions when mixture models are utilized for the parameter densities of the system. In addition, it is demonstrated how the stability analysis based on these randomly distributed parameters can be performed. Several interpretations of such systems and the corresponding results are provided, including the concept of quantum-like methods in relation to incoherent superposition, which we regard as a legitimate and novel way of looking at such systems. It is demonstrated that all numerical calculations can be performed efficiently utilizing recently developed methods. Future work may include the application of this probabilistic steady-state analysis framework to other highly nonlinear dynamical systems, investigating their uncertainty propagation and superposition of parameters.

\subsection*{Declaration of generative AI use}

Generative AI (OpenAI, GPT-5.6 Terra) was used to improve the spelling and grammar of the manuscript.

% Bibliography

\end{document}